\input amstex
\documentstyle{amsppt}

\topmatter

\title  Philip Hall's 'rather curious' formula and the Rogers-Ramanujan identities  \endtitle

\author Avinoam Mann\endauthor

\endtopmatter

\bigskip

In {\bf [H1]}, Philip Hall published the following remarkable identity, which, in a typical understatement, he described as "rather curious". Fix a
prime $p$, and let $\Cal A_p$ be the class of all finite abelian $p$-groups. Then

$$ \sum_{G\in\Cal A_p} \frac{1}{|G|} = \sum_{G\in\Cal A_p} \frac{1}{|Aut(G)|}.\leqno {(1)}$$

Since the number of abelian groups of order $p^n$ is $\pi (n)$, the number of partitions of $n$, another way of
writing (1) is

$$\sum_G \frac{1}{|Aut(G)|} = \sum \frac{\pi(n)}{p^n}.\leqno {(2)}$$

Alternative proofs to Hall's were given by I.G.Macdonald {\bf [Mc]} (1984, applying {\it Hall algebras}),
H.Cohen-H.W.Lenstra {\bf [CL]} (1984, with a number-theoretical motivation), T.Yoshida {\bf [Y]} (1992), the
author {\bf [M1]} (1996), F.Clarke {\bf [C]} (2006), and J.Lengler {\bf [L]} (2008). Also, since (1) can be expressed as a partition identity, or a so-called q-identity, it can be derived from other identities of that type (see e.g. {\bf {[St]}}), but in papers doing so the group theoretical connection is often not made explicit \footnote {I am grateful to O.Warnaar for pointing me to these references.}. In this note we point out how (1) can be considered as the limiting form of a generalized Rogers-Ramanujan identity.

Recall that the Rogers-Ramanujan identities are

$$1+\sum_{k=1}^{\infty}\frac{x^{k^2}}{(1-x)(1-x^2)...(1-x^k)}=\Pi_{j=0}^{\infty}\frac{1}{(1-x^{5j+1})(1-x^{5j+4})},\leqno{(3)}$$

$$1+\sum_{k=1}^{\infty}\frac{x^{k^2+k}}{(1-x)(1-x^2)...(1-x^k)}=\Pi_{j=0}^{\infty}\frac{1}{(1-x^{5j+2})(1-x^{5j+3})}.\leqno{(4)}$$

Here $x$ can be any number such that $|x| < 1$ (or we can consider the identities as formal ones).

There seem to be two mysteries in the R-R identities. First, why should an infinite series be equal to an infinite product? and second, why the special role of the number 5? This latter mystery is somewhat dispelled by a generalization, which shows that the R-R identities are just the simplest instance of an infinite series of identities, where 5 is replaced by all odd numbers, taken in turn. These are

$$\sum \frac {x^{k_1^2+...+k_r^2+k_i+...+k_r}}{(1-x)(1-x^2)...(1-x^{k_1-k_2})...(1-x)...(1-x^{k_{r-1}-k_r})(1-x)...(1-x^{k_r})}\leqno{(5)}$$

$$ = \prod_{j>0,~j\not\equiv 0,\pm i~(mod~2r+1)}\frac{1}{(1-x^j)}.$$

\noindent Here $r$ and $i$ are fixed, $r \ge 1$, $1 \le i \le r+1$, $k_1 \ge k_2 \ge ... \ge k_r \ge 0$ {\bf [An, 7.8]}. The R-R identities are obtained upon taking $r = 1$ and $i = 2$ or $i = 1$.

To see the connection to Hall's formula, write an abelian group $G$ of order $p^n$ as a direct sum of cyclic subgroups of orders $p^{\lambda_1}, ..., p^{\lambda_r}$,
where $\lambda_1 \ge ... \ge \lambda_r$ and $n = \sum~ \lambda_i$. To the partition $(\lambda_1, ..., \lambda_r)$ we associate its Young diagram, which consists of $r$ rows, arranged with the longest one on top, and the $i$th one consists of $\lambda_i$ unit squares. The dual partition $\mu_1, ..., \mu_s$ is obtained by counting the number of squares in columns, not rows. Then

$$ |Aut(G)| = \frac{f_{\mu_1-\mu_2}(\rho)f_{\mu_2-\mu_3}(\rho)...}{\rho^{\mu_1^2+\mu_2^2+ ...}},\leqno{(6)}$$

where $\rho = 1/p$ and $f_k(x) = (1-x)(1-x^2) ... (1-x^k)$. \footnote {For this formula, Hall refers to {\bf [Sp]}. However, the elegant and suggestive form (6) is due to him.} If we take in (5) $i = k+1$ and $x = 1/p$, the general term on the LHS is $\frac{1}{|Aut(G)|}$, where $G$ is the abelian $p$-group determined by the dual partition to $(k_1,...,k_r)$, and if we let $r$ tend to infinity, the RHS becomes the generating function $\sum \pi(n)x^n$ of all partitions, and (5) becomes (1). Thus the R-R identities and Hall's formula are the two extremes, $r = 1$ and $r = \infty$, of (5). If we keep $r$ fixed in (5), the LHS becomes $\sum_{G\in \Cal A_p,~exp(G) \le p^r} \frac{1}{|Aut(G)|}$, a fact that was noticed in {\bf [L, 7.1]}. If we take $i = 1$, the general term on the LHS becomes $\frac{1}{|G||Aut(G)|} = \frac{1}{|Hol(G)|}$, where $G$ is as above, and $Hol(G)$ is the {\it holomorph} of $G$, its semidirect product with $Aut(G)$ (acting naturally). Letting again $r$ tend to infinity, we get on the RHS $\prod_{j \ge 2} \frac{1}{1-x^j}$, which is the generating function for the number of partitions with components at least 2. By duality, the number of such partitions equals the number of partitions with two equal largest components, which partitions correspond to {\it capable} abelian $p$-groups, i.e. groups that are isomorphic to the central factor group $H/Z(H)$ of some group. Thus, letting $\Cal B_p$ be the class of capable abelian $p$-groups, we obtain, this time as the analogue of the second R-R formula

$$ \sum_{G\in\Cal B_p} \frac{1}{|G|} =  \sum_{G\in\Cal A_p}\frac{1}{|Hol(G)|}.\leqno {(7)}$$

In {\bf [M3]} we show that the sum in (7) equals $\frac{p-1}{p} \sum_{G\in\Cal A_p} \frac{1}{|G|}$. This, and (7), occur as the special case $k = 1$ of: for $k \ge 0$,

$$\sum_{G\in\Cal A_{p,k+1}}\frac{1}{|G|} = \sum_{G\in\Cal A_p}\frac{1}{|G|^k|Aut(G)|}
= (1-\frac{1}{p})(1-\frac{1}{p^2})...(1-\frac{1}{p^k})\sum_{G\in\Cal A_p} \frac{1}{|G|}.\leqno{(8)}$$

\noindent  Here $\Cal A_{p,k}$ is the class of abelian $p$-groups, in which, in the above decomposition, the first $k$ factors are isomorphic. The case $k = 0$ of (8) is Hall's equation (1). The proof of (8) is by a variation of Hall's original argument.

\smallskip

{\bf Notes. 1.} There is a 'non-elegant' variation of (8) involving terms of the form $\frac{|G|^k}{|Aut(G)|}$.

{\bf 2.} It is shown in Proposition 17 of {\bf [M2]} that $\sum\frac {\pi(n)}{p^n}$ is irrational.

\bigskip

\head References \endhead

{\bf An.} G.E.Andrews, {\it The theory of partitions}, Addison-Wesley, Reading (Mass.), 1976.

{\bf C.} F.W.Clarke, Counting abelian group structures, Proc. Amer. Math. Soc. 134 (2006), 2795-2799.

{\bf CL.} H.Cohen-H.W.Lenstra, Heuristics on class groups of number fields, in {\it Number Theory, Noordwijkerhout
1983}, Lecture Notes in Mathematics 1068, Springer, Berlin 1984, 33-62.

{\bf H1.} P.Hall, A partition formula connected with abelian groups, Comm. Math. Helv. 11 (1938), 126-129 (= {\bf
[H4]}, 245-248).

{\bf H2.} K.W.Gruenberg-J.E.Roseblade (editors), {\it Collected Works of Philip Hall}, Clarendon Press, Oxford
1988.

{\bf HW.} G.H.Hardy-E.M.Wright, {\it An Introduction to the Theory of Numbers}, 56th ed., Oxford University Press,
Oxford 2008.

{\bf L.} J.Lengler, A combinatorial interpretation of the probabilities of $p$-groups in the Cohen-Lenstra
measure, J. Number Theory 128 (2008), 2070-2084.

{\bf Mc.} I.G.Macdonald, The algebra of partitions, in {\it Group Theory: Essays for Philip Hall}, Academic Press,
London 1984, 315-333.

{\bf M1.} A.Mann, Philip Hall's `rather curious' formula for abelian $p$-groups, Israel J. Math. (S.A.Amitsur
Memorial Volume) 96 (1996), 445-448.

{\bf M2.} A.Mann, Subgroup growth in pro-p groups, in {\it New Horizons in pro-p Groups}, edited by M.P.F.du
Sautoy, D.Segal, and Aner Shalev, Birkhauser, Boston 2000, 233-247.

{\bf M3.} A.Mann, Some group theoretical mass formulae, in preparation.

{\bf Sp.} A.Speiser, {\it Die Theorie der Gruppen von endlicher Ordnung}, 3rd ed., Dover, New York 1943.

{\bf St.} J.R.Stembridege, Hall-Littlewood functions, plane partitions, and the Rogers-Ramanujan identities. Trans. Amer. Math. Soc.  319  (1990), 469–498.

{\bf Y.} T.Yoshida, P.Hall's strange formula for abelian $p$-groups, Osaka Math. J. 29 (1992), 421-431.

\end